\documentclass{article}
\usepackage{eurosym}
\usepackage{amssymb}
\usepackage{amsfonts}
\usepackage{amsmath}
\usepackage{mathrsfs}

\newtheorem{theorem}{Theorem}

\newtheorem{corollary}[theorem]{Corollary}

\newtheorem{definition}[theorem]{Definition}

\newtheorem{lemma}[theorem]{Lemma}

\newtheorem{remark}[theorem]{Remark}

\newenvironment{proof}[1][Proof]{\noindent\textbf{#1.} }{\ \rule{0.5em}{0.5em}}

\begin{document}
\title{\textbf{Riemannian Submersions Whose Total Manifolds Admitting a Ricci Soliton}}
\author{\c{S}emsi Eken Meri\c{c} and Erol K{\i}l{\i}\c{c}  }
\date{\vspace{-5ex}}
\maketitle

\bigskip
\begin{center}
	$^{1}$ Department of Mathematics, Mersin University, 33343, Mersin, Turkey  \\
		$^{2}$ Department of Mathematics, \.{I}n\"{o}n\"{u} University, 44280, Malatya, Turkey\\
	 $^{1}$semsieken@hotmail.com, $^{2}$erol.kilic@inonu.edu.tr
\end{center}

%\date{}  % Toggle commenting to test

\bigskip

\textbf{Abstract} In this paper, we study Riemannian submersions whose total manifolds admitting a Ricci soliton. Here, we characterize any fiber of such a submersion is Ricci soliton or almost Ricci soliton. Indeed, we obtain necessary conditions for which the target manifold of Riemannian submersion is a Ricci soliton. Moreover, we study the harmonicity of Riemannian submersion from Ricci soliton and give a characterization for such a submersion to be harmonic.

\bigskip

 \medskip


\medskip

\section{Introduction}

\noindent In \cite{Hamilton}, Hamilton defined the notion of Ricci flow and showed that the self similar solutions of such a flow are Ricci solitons. According to the definition of Hamilton, a Riemannian manifold $(M,g)$ is said to be a Ricci soliton if it satisfies
\begin{equation}
\frac{1}{2}\mathscr{L}_Vg+Ric+\lambda g=0  ,  \label{p:1}
\end{equation}
where $\mathscr{L}_Vg$ is the Lie-derivative of the metric tensor of $g$ with respect to $V$, $Ric$ is the Ricci tensor of $(M,g)$, $V$ is a vector field (the potential field) and $\lambda$ is a constant on $M$. We shall denote a Ricci soliton by $(M,g,V,\lambda)$. The Ricci soliton $(M,g,V,\lambda)$ is said to be shrinking, steady or expanding acoording as $\lambda<0$, $\lambda=0$ or $\lambda>0$, respectively.\medskip

The Ricci soliton $(M,g,V,\lambda)$ is said to be a graident Ricci soliton, if the potential field $V$ is the gradient of some smooth function $f$ on $M$, which is denoted by $(M,g,f,\lambda)$. On the other hand, a trivial Ricci soliton is an Einstein metric if the potential field $V$ is zero or Killing, that is, the Lie-derivative of the metric tensor $g$ with respect to $V$ vanishes, identically. Moreover, a tangent vector field $V$ on a Riemannian manifold $(M,g)$ is called conformal, if it satisfies
\begin{eqnarray*}
	\mathscr{L}_{V}g & = & 2fg,
\end{eqnarray*}
where $\mathscr{L}_Vg$ is the Lie-derivative of $(M,g)$ and $f$ is called the potential function of $V$.\medskip

Pigola et al. defined a new class of Ricci solitons by taking $\lambda$ is a variable function instead of the constant and then, the Ricci soliton $(M,g,V,\lambda)$ is called an almost Ricci soliton. Hence, the almost Ricci soliton becomes a Ricci soliton, if the function $\lambda$ is a constant. Indeed, the almost Ricci soliton is called shrinking, steady or expanding acoording as $\lambda<0$, $\lambda=0$ or $\lambda>0$, respectively (\cite{Pigola}).\medskip

In \cite{Perelman}, Perelman used the Ricci soliton in order to solve the Poincar\'{e} conjecture and then the geometry of Ricci solitons has been the focus of attention of many mathematicians. For example, Chen and et al. studied the Riemannian manifolds endowed with a concurrent, concircular or torqued vector fields admitting a Ricci soliton and gave many characterizations [4-8].\medskip

Moreover, Ricci solitons have been studied on contact, paracontact and Sasakian manifolds. For instance, Blaga and et al. focused on $\eta-$Ricci solitons and provided some remarks in Sasakian manifolds (\cite{Blaga:1}, \cite{Blaga:2}). Also, in \cite{Selcen}, Perkta\c{s} and et al. gave some characterizations on Ricci solitons in 3-dimensional normal almost para-contact metric manifolds.\medskip

On the other hand, the theory of Riemannian submersions is a very interesting topic in differential geometry, 
since it has many applications in physics and mechanics. For instance, in the theory of 
Kaluza-Klein, one starts with the hypothesis that the space-time has \textit{4+m}-dimensions. 
An interesting mechanism for space-time compactification is proposed in the form of a 
non-linear sigma model. The general solutions of this model can be expressed in terms of 
harmonic maps satisfying the Einstein equations. A very general class of solutions is given 
by Riemannian submersions of this model (see \cite{r6}, chapter 9). \medskip

The concept of Riemannian submersion between Riemannian manifolds was initated by O'Neill and Gray, independently. They gave some basis formulas on the theory of Riemannian submersions and it has been extended in the last three decades. (For the theory of Riemannian submersions, we refer \cite{r7}, \cite{r9} and \cite{sahin}). \medskip

In recent years, Meri\c{c} and et al. established many inequalities for Riemannian submersions to obtain the relationships between the intrinsic and extrinsic invariants for such a submersion (\cite{Eken:1}, \cite{Eken:3} and \cite{Gulbahar}).\medskip

This work is organized as follows: The section 2 is a very brief of review of Riemannian submersion. In Section 3, a Riemannian submersion $\pi$ from Ricci soliton to a Riemannian manifold is considered and the necessary conditions for which any fiber of $\pi$ is a Ricci soliton or an almost Ricci soliton are given. The last section is devoted to harmonicity. Using the tension field, a necessary and sufficient condition for which a Riemannian submersion from Ricci soliton is harmonic is obtained.

\section{Preliminaries}

\noindent In this section, we recall some basic notions about Riemannian submersions from \cite{r6}: \medspace

A map $\pi:(M,g)\rightarrow (B,g^{'}) $ is called a $C^{\infty}$-submersion between Riemannian manifolds $(M^m,g)$ and $(B^n,g^{'})$, if $\pi$ has a maximal rank at any point of $M$. For any $x\in B$, $\pi^{-1}(x)$ is closed r-dimensional submanifold of $M$, such that 
$r=m-n$. For any $p\in M$, denoting $\mathcal{V}_{p}=ker{\pi_{*p}}$ and then the distribution $\mathcal{V}$ is integrable. Also, $T_{p}\pi^{-1}(x)$ are $r-$dimensional subspaces of $\mathcal{V}_{p}$ and it follows that $\mathcal{V}_{p}=%
T_{p}\pi^{-1}(x)$. Hence, $\mathcal{V}_{p}$ is called the vertical space of any point
$p\in M$. \\
\noindent Denote the complementary distribution of $\mathcal{V}$ by $\mathcal{H}$, the one has 
\begin{eqnarray*}
	T_{p}(M)=\mathcal{V}_{p}\oplus \mathcal{H}_{p} \label{crv:q}
\end{eqnarray*}
where $\mathcal{H}_{p}$ is called the horizontal space of any point $p\in M$.
\smallskip

Let $ \pi:(M,g)\rightarrow(B,g^{'})$ be a submersion between Riemannian manifolds $(M,g)$ and $(B,g^{'})$. At any point $p\in M$, we say that $\pi$ is a Riemannian submersion if $\pi_{*p}$ preserves the length of the horizontal vectors.\medskip

\noindent Some basic properties about Riemannian submersion are presented as follows:\medskip

Let $\pi:(M,g)\rightarrow(B,g^{'})$ be a Riemannian submersion, 
and denote by $\nabla$ and $\nabla^{'}$ the Levi-Civita connections of $M$ and $B$, 
respectively. If $X$, $Y$ are the basic vector fields, $\pi$-related to $X^{'}, Y^{'},$ 
one has:\medskip

\noindent		$ (i) \ \ g(X,Y)=g^{'}(X^{'},Y^{'})\circ \pi,$\\
$ (ii) \ \ h[X,Y]$ is the basic vector field $\pi$-related to $[X^{'},Y^{'}]$,\\
$ (iii)\ h(\nabla_{X}Y)$ is the basic vector field $\pi$-related to 
$\nabla^{'}_{X^{'}}Y^{'}$,\\
$(iv)$ for any vertical vector field $V$, $[X,V]$ is the vertical. \medskip

A Riemannian submersion $\pi:(M,g)\rightarrow(B,g^{'})$ determines two tensor fields $\mathcal{T}$ and $\mathcal{A}$ on the base manifold $M$, which are called the fundamental tensor fields or the invariants of Riemannian submersion $\pi$ and they are defined by
\begin{eqnarray*}
	&&\mathcal{T}(E,F)=\mathcal{T}_{E}F=h\nabla_{vE}vF+v\nabla_{vE}hF, \\
	&&\mathcal{A}(E,F)=\mathcal{A}_{E}F=v\nabla_{hE}hF+h\nabla_{hE}vF,
\end{eqnarray*}
where $v$ and $h$ are the vertical and horizontal projections, respectively and $\nabla$ is a Levi-Civita connection of $M$, for any $E,F\in \Gamma(TM)$. Indeed, the fundamental tensors $\mathcal{T}$ and $\mathcal{A}$ satisfy the followings:\medskip

\noindent$ (i) \ \ \mathcal{T}_{V}W=\mathcal{T}_{W}V,$ \\
$ (ii) \ \ \mathcal{A}_{X}Y=-\mathcal{A}_{Y}X=\frac{1}{2}v[X,Y]$,   \medskip

\noindent for any $V,W \in \Gamma(VM)$ and $X,Y \in \Gamma(HM)$.\medskip

We here note that the vanishing of the tensor field $\mathcal{A}$ means the horizontal distribution $\mathcal{H}$ is integrable. On the other hand, the vanishing of the tensor field $\mathcal{T}$ means any fibre of Riemannian submersion $\pi$ is totally geodesic submanifold of $M$. Also, any fiber of Riemannian submersion $\pi$ is totally umbilical if and only if
\begin{equation*}
\mathcal{T}_VW=g(V,W)H,
\end{equation*}
where $H$ denotes the mean curvature vector field of any fiber in $M$, for any $V,W\in\Gamma(VM)$.\medskip

\noindent Moreover, for any $E,F,G\in\Gamma(TM)$ one has
\begin{eqnarray}
g(\mathcal{T}_EF,G)+g(\mathcal{T}_EG,F)=0,  \label{u}\ \ \\ g(\mathcal{A}_EF,G)+g(\mathcal{A}_EG,F)=0. \label{key:43}
\end{eqnarray}

\noindent Using fundamental tensor fields $\mathcal{T}$ and $\mathcal{A}$, the following formulas are given as
\begin{eqnarray}
\nabla_{V}W=\mathcal{T}_{V}W+\hat{\nabla}_{V}W, \label{re:1}\\
\nabla_{V}X=h\nabla_{V}X+\mathcal{T}_{V}X, \label{re:2}\\
\nabla_{X}V=\mathcal{A}_{X}V+v\nabla_{X}V, \label{re:3}\\
\nabla_{X}Y=h\nabla_{X}Y+\mathcal{A}_{X}Y \label{re:4}
\end{eqnarray}
\textit{for any} $V,W\in \Gamma(VM)$ \textit{and} $X,Y\in \Gamma(HM)$. 
\medskip

\noindent Denoting the Riemannian curvature tensors of $(M,g)$, $(B,g^{'})$ and any fibre of $\pi$ by $R$, $R^{'}$ and  $\hat{R}$ respectively. Then, 
\begin{eqnarray*}
	R(U,V,F,W)&=& \hat{R}(U,V,F,W)-g(\mathcal{T}_{U}W,\mathcal{T}_{V}F)+g(\mathcal{T}_{V}W,\mathcal{T}_{U}F),\label{crv:3} \\
	R(X,Y,Z,H)&=& R^{'}(X^{'},Y^{'},Z^{'},H^{'})\circ \pi+2g(\mathcal{A}_{X}Y,\mathcal{A}_{Z}H) \label{crv:35}\\
	&&-g(\mathcal{A}_{Y}Z,\mathcal{A}_{X}H)+g(\mathcal{A}_{X}Z,\mathcal{A}_{Y}H), \nonumber 
\end{eqnarray*}
for any $U,V,W,F\in \Gamma(VM)$ and $X,Y,Z,H\in \Gamma(HM)$. 
\medskip

\noindent Let $\{U,V\}$ and $\{X,Y\}$ be an orthonormal basis of the vertical and horizontal $2-$plane, respectively. Then, one has

\begin{eqnarray*}
	K(U,V)& = & \hat{K}(U,V)-\|\mathcal{T}_{U}V\|^{2}+g(\mathcal{T}_{U}U,\mathcal{T}_{V}V),\\
	K(X,Y) & = & K^{'}({\pi_{*}X,\pi_{*}Y})+3\|\mathcal{A}_{X}Y\|^{2},
\end{eqnarray*}
where $K$, $\hat{K}$ and $K^{'}$ are the sectional curvature in the total space $M$, any fiber of $\pi$ and $(B,g^{'})$, respectively.\\

The Ricci tensor $Ric$ on $(M,g)$ is given by 

\begin{eqnarray}
Ric(U,W) & = & \hat{Ric}(U,W)+g(N,\mathcal{T}_{U}W)-\sum_{i=1}^ng((\nabla_{X_{i}}\mathcal{T})(U,W),{X_{i}})   \label{sd:1} \\
&& -g(\mathcal{A}_{X_{i}}U,\mathcal{A}_{X_{i}}W) \notag \\
Ric(X,Y) & = & Ric^{'}(X^{'},Y^{'})\circ \pi -\frac{1}{2}\{g(\nabla_{X}N,Y)+g(\nabla_{Y}N,X)\}  \label{sd:2} \\
&& +2\sum_{i=1}^{n}g(\mathcal{A}_{X}X_{i},\mathcal{A}_{Y}X_{i})+\sum_{j=1}^{r}g(\mathcal{T}_{U_{j}}X,\mathcal{T}_{U_{j}}Y)	 \notag \\
Ric(U,X) & = & -g(\nabla_{U}N,X)+\sum_{j}g((\nabla_{U_{j}}\mathcal{T})(U_{j},U),X) \label{sd:3} \notag \\
&& -\sum_{i=1}^{n}\{ g((\nabla_{X_{i}}\mathcal{A})(X_{i},X),U)+2g(\mathcal{A}_{X_{i}}X,\mathcal{T}_{U}X_{i})\} \notag
\end{eqnarray}
where $\{X_i\}$ and $\{U_j\}$ are the orthonormal basis of $\mathcal{H}$ and $\mathcal{V}$, respectively for any $U,V\in\Gamma(VM)$ and $X,Y\in \Gamma(HM)$.\\

\noindent On the other hand, the mean curvature vector field $H$ on any fibre 
of Riemannian submersion $\pi$ is given by
\begin{eqnarray*}
	N=rH,
\end{eqnarray*}
such that
\begin{eqnarray}
N=\sum_{j=1}^{r}\mathcal{T}_{U_{j}}U_{j} \label{bbb:1}
\end{eqnarray}
and $r$ denotes the dimension of any fibre of $\pi$ and $\{U_{1},U_{2},...,U_{r}\}$ 
is an orthonormal basis on vertical distribution. We remark that 
the horizontal vector field $N$ vanishes if and only 
if any fibre of Riemannian submersion $\pi$ is minimal.\\

\noindent Furthermore, using the equality (\ref{bbb:1}), we get
\begin{eqnarray*}
	g(\nabla_{E}N,X)=\sum_{j=1}^{r}g((\nabla_{E}\mathcal{T})(U_{j},U_{j}),X) \label{j:1}
\end{eqnarray*}
for any $E\in \Gamma(TM)$ and $X\in \Gamma(HM)$.
\medskip

\noindent We denote the horizontal divergence of any vector field $X$ on 
$\Gamma(HM)$ by $\check{\delta}(X)$ and given by 
\begin{eqnarray*}
	\check{\delta}(X) & = & \sum_{i=1}^{n}g(\nabla_{X_{i}}X,X_{i}), \label{bbb:2}
\end{eqnarray*}
where $\{X_{1},X_{2},...,X_{n}\}$ is an orthonormal basis of horizontal space 
$\Gamma(HM)$. Hence, considering (\ref{bbb:1}), we have
\begin{eqnarray}
\check{\delta}(N) & = & \sum_{i=1}^{n}\sum_{j=1}^{r}%
g((\nabla_{X_{i}}\mathcal{T})(U_{j},U_{j}),X_{i}). \label{bb:3}
\end{eqnarray}
\noindent For details, we refer to (\cite{Bessel}, pp. 243).

\section{Riemannian Submersions Whose Total Manifolds Admitting a Ricci Soliton}

In the present section, we study a Riemannian submersion $\pi:(M,g)\rightarrow (B,g')$ from the Ricci soliton onto a Riemannian manifold and give some characterizations for any fiber of such a submersion and the target manifold $B$.\\

\noindent Using Eqs. (2)-(7) for Riemannian submersions, we give the following lemma:

\begin{lemma} \label{l:1}
	Let $\pi:(M,g)\rightarrow (B,g')$ be a Riemannian submersion between Riemannian manifolds. Then the followings are equaivalent to each other:\\
	
	\noindent(i) the vertical distribution $\mathcal{V}$ is parallel,\\
	(ii) the horizontal distribution $\mathcal{H}$ is parallel,\\
	(iii) the fundamental tensor fields $\mathcal{T}$ and $\mathcal{A}$ are vanish identically, that is, $\mathcal{T}\equiv 0$ and $\mathcal{A}\equiv 0$.
\end{lemma}

\begin{theorem}	
	Let $(M,g,V,\lambda)$ be a Ricci soliton with the vertical potential field $V$ and $\pi:(M,g)\rightarrow (B,g')$ be a Riemannian submersion between Riemannian manifolds. If one of the conditions $(i)-(iii)$ in Lemma \ref{l:1} is satisfied, then any fibre of Riemannian submersion $\pi$ is a Ricci soliton.
\end{theorem}

\begin{proof}  
	Since $(M,g)$ is a Ricci soliton, one has
	\begin{eqnarray*}
		\frac{1}{2}(\mathscr{L}_{V}g)(U,W)+Ric(U,W)+\lambda g(U,W)=0,
	\end{eqnarray*} 
	for any $U,W\in\Gamma(VM)$. Using the equality (\ref{sd:1}), we have
	\begin{eqnarray}
	&& \frac{1}{2} \big\{g(\nabla_{U}V,W)+g(\nabla_{W}V,U)  \big\}+\hat{Ric}(U,W)+g(N,\mathcal{T}_{U}W) \label{sd:4} \\
	&& -\sum_{i=1}^{n}\big( g((\nabla_{X_{i}}\mathcal{T})(U,W),X_{i})-g(\mathcal{A}_{X_{i}}U,\mathcal{A}_{X_{i}}W)\big) +\lambda g(U,W)=0, \notag
	\end{eqnarray}
	where $\{X_i\}$ denotes an orthonormal basis of the horizontal distribution $\mathcal{H}$ and $\nabla$ is the Levi-Civita connection on $M$. Using Lemma \ref{l:1} and the equality (\ref{re:1}), it follows
	\begin{eqnarray*}
		\frac{1}{2} \{{g}(\hat{\nabla}_{U}V,W)+{g}(\hat{\nabla}_{W}V,U)  \}+\hat{Ric}(U,W)+\lambda {g}(U,W)=0,
	\end{eqnarray*}
	which means any fibre of the Riemannian submersion $\pi$ is a Ricci soliton.
\end{proof}

\begin{theorem}	\label{l:3}
	Let $(M,g,V,\lambda)$ be a Ricci soliton with the vertical potential field $V$ and $\pi:(M,g)\rightarrow (B,g')$ be a Riemannian submersion between Riemannian manifolds with totally umbilical fibres. If the horizontal distribution $\mathcal{H}$ is integrable, then any fibre of Riemannian submersion $\pi$ is an almost Ricci soliton.
\end{theorem}

\begin{proof}
	Since the total space $(M,g)$ of Riemannian submersion $\pi$ admits a Ricci soliton, then using (\ref{p:1}) and (\ref{sd:1}), we have
	\begin{eqnarray}
	&& \frac{1}{2} \big\{g(\nabla_{U}V,W)+g(\nabla_{W}V,U)  \big\}+\hat{Ric}(U,W)+\sum_{j=1}^{r}g(\mathcal{T}_{U_j}U_j,\mathcal{T}_{U}W)  \label{bbb:3}\\
	&&-\sum_{i=1}^{n}(g((\nabla_{X_{i}}\mathcal{T})(U,W),X_i)-g(\mathcal{A}_{X_{i}}U,\mathcal{A}_{X_{i}}W))+\lambda g(U,W)=0,  \notag
	\end{eqnarray}
	for any $U,W\in\Gamma(VM)$. Also, the Ricci soliton $(M,g,V,\lambda)$ has totally umbilical fibres and putting (\ref{re:1}) in (\ref{bbb:3}), one has
	\begin{eqnarray}
	&& \frac{1}{2} \{g(\hat{\nabla}_{U}V,W)+g(\hat{\nabla}_{W}V,U)  \}+\hat{Ric}(U,W)+\sum_{j=1}^rg(\mathcal{T}_{U_j}U_j,\mathcal{T}_{U}W) \notag\\
	&&-\sum_{i=1}^{n}\{(\nabla_{X_{i}}g)(U,W)g(H,X_i)-g(\nabla_{X_{i}}H,X_i){g}(U,W)	\} \notag \\
	&& -\sum_{i=1}^{n}g(\mathcal{A}_{X_{i}}U,\mathcal{A}_{X_{i}}W)+\lambda g(U,W)=0.	\notag
	\end{eqnarray} 	
	Since the horizontal distribution $\mathcal{H}$ is integrable, we have
	\begin{eqnarray*}
		&& \frac{1}{2}(\mathscr{L}_V{g})(U,W)+\hat{Ric}(U,W)-\sum_{i=1}^ng(\nabla_{X_{i}}H,X_i){g}(U,W)+r\Vert H  \Vert^2{g}(U,W)\\
		&&+\lambda {g}(U,W)=0,
	\end{eqnarray*}
	where $H$ is the mean curvature vector of any fibre of $\pi$. From (\ref{bbb:2}), we obtain	
	\begin{eqnarray}
	\frac{1}{2}(\mathscr{L}_V{g})(U,W)+\hat{Ric}(U,W)+ \big( r\Vert H \Vert^2-\check{\delta}(H)+\lambda \big){g}(U,W) & = & 0,
	\end{eqnarray}
	which means any fibre of $\pi$ is an almost Ricci soliton.	
\end{proof}

Considering Theorem \ref{l:3}, we get the following:
\begin{corollary}	
	Let  $(M,g,V,\lambda)$ be a Ricci soliton and $\pi:(M,g)\rightarrow (B,g^{'})$ be a Riemannian submersion between Riemannian manifolds, such that the horizontal distribution $\mathcal{H}$ is integrable. Any fiber of $\pi$ is a Ricci soliton, if one of the following conditions satisfies: \medskip
	
	\noindent  (i) Any fiber of $\pi$ is a totally umbilical and has a constant mean curvature.\\
	\noindent (ii) Any fiber of $\pi$ is a totally geodesic.
\end{corollary}

Then, we have the following theorem:

\begin{theorem} \label{tt:1}
	Let $(M,g,E,\lambda)$ be a Ricci soliton with the potential field $E\in\Gamma(TM)$ and $\pi:(M,g)\rightarrow (B,g')$ be a Riemannian submersion between Riemannian manifolds. If one of the conditions $(i)-(iii)$ in Lemma \ref{l:1} is satisfied, then the followings are satisfied:\medskip
	
	\noindent (i) If the vector field $E$ is vertical, then  $(B,g^{'})$ is an Einstein.\\	
	\noindent (ii) If the vector field $E$ is horizontal, then  $(B,g^{'})$ is a Ricci soliton with potential field $E^{'}$, such that $\pi_{*}E=E^{'}$. 
\end{theorem}

\begin{proof} Since the total space $(M,g)$ of Riemannian submersion $\pi$ admits a Ricci soliton with potential field $E\in\Gamma(TM)$, then using (\ref{p:1}) and (\ref{sd:2}), we have
	
	\begin{eqnarray}
	&& \frac{1}{2} \{g(\nabla_{X}E,Y)+g(\nabla_{Y}E,X)  \}+Ric^{'}(X^{'},Y^{'})\circ \pi
	-\frac{1}{2}\big( g(\nabla_{X}N,Y) \label{sd:5} \\ 
	&& +g(\nabla_{Y}N,X)\big) 
	+2\sum_{i=1}^ng(\mathcal{A}_{X}X_{i},\mathcal{A}_{Y}X_{i})+\sum_{j=1}^r g(\mathcal{T}_{U_{j}}X,\mathcal{T}_{U_{j}}Y) \notag \\
	&& +\lambda g(X,Y)=0, \notag
	\end{eqnarray}
	where $X^{'}$ and $Y^{'}$ are $\pi-$related to $X$ and $Y$ respectively, for any $X,Y\in\Gamma(HM)$.\\	
	
	\noindent	Applying Lemma \ref{l:1} to above equation (\ref{sd:5}), one has	
	\begin{eqnarray}
	&&\frac{1}{2} \big\{ g(\nabla_{X}E,Y)+g(\nabla_{Y}E,X)  \big\}+Ric^{'}(X^{'},Y^{'})\circ \pi +\lambda g(X,Y)=0.       \label{l:25}
	\end{eqnarray}	
	
	\noindent	(i) If the vector field $E$ is vertical and using the equality (\ref{re:3}) in (\ref{l:25}), we get
	\begin{eqnarray*}
		&&\frac{1}{2} \big\{ g(\mathcal{A}_{X}E,Y)+g(\mathcal{A}_{Y}E,X)  \big\}+Ric^{'}(X^{'},Y^{'})\circ \pi +\lambda g(X,Y)=0.       
	\end{eqnarray*}
	Since one of the conditions $(i)-(iii)$ in Lemma \ref{l:1} is satisfied, one has
	\begin{equation*}
	Ric^{'}(X^{'},Y^{'})\circ \pi +\lambda g(X,Y)=0.
	\end{equation*}
	Because $X$ and $Y$ are $\pi-$related to $X^{'}$ and $Y^{'}$, we have
	\begin{eqnarray*}
		Ric^{'}(X^{'},Y^{'})\circ \pi +\lambda g^{'}(X^{'},Y^{'})\circ \pi=0,
	\end{eqnarray*}
	which is equivalent to
	\begin{eqnarray*}
		Ric^{'}(X^{'},Y^{'})+\lambda g^{'}(X^{'},Y^{'})=0.
	\end{eqnarray*}
	Therefore $(B,g^{'})$ is an Einstein.\\
	
	\noindent (ii) If the vector field $E$ is a horizontal, the equation (\ref{l:25}) follows
	\begin{eqnarray}
	\frac{1}{2}(\mathscr{L}_Eg)(X,Y)+Ric^{'}(X^{'},Y^{'})\circ \pi+\lambda g(X,Y)=0. \label{key:s}
	\end{eqnarray}
	On the other hand, $h(\nabla_{X}E)$ and $h(\nabla_{Y}E)$ are the basic vector fields $\pi-$related to $\nabla^{'}_{X^{'}}E^{'}$ and $\nabla^{'}_{Y^{'}}E^{'}$ respectively, the equation (\ref{key:s}) gives
	\begin{eqnarray*}
		\frac{1}{2}(\mathscr{L}_{E^{'}}g^{'})(X^{'},Y^{'})\circ \pi+Ric^{'}(X^{'},Y^{'})\circ \pi+\lambda g^{'}(X^{'},Y^{'})\circ \pi=0,
	\end{eqnarray*}
	which means
	\begin{eqnarray*}
		\frac{1}{2}(\mathscr{L}_{E^{'}}g^{'})(X^{'},Y^{'})+Ric^{'}(X^{'},Y^{'})+\lambda g^{'}(X^{'},Y^{'})=0.
	\end{eqnarray*}
	Therefore the proof is completed.
\end{proof}

Using Lemma \ref{l:1} and the equality (\ref{sd:2}), we get the following:

\begin{remark}
	Let $(M,g,\xi,\lambda)$ be a Ricci soliton with the horizontal potential field $\xi$ and $\pi:(M,g)\rightarrow (B,g')$ be a Riemannian submersion between Riemannian manifolds. If one of the conditions $(i)-(iii)$ in Lemma \ref{l:1} is satisfied, the the vector field $N$ is Killing on the horizontal distribution $\mathcal{H}$.
\end{remark}

\begin{theorem} \label{prop:12}
	Let $(M,g,\xi,\lambda)$ be a Ricci soliton with the horizontal potential field $\xi$ and $\pi:(M,g)\rightarrow (B,g')$ be a Riemannian submersion from Riemannian manifold to an Einstein manifold. If one of the conditions $(i)-(iii)$ in Lemma \ref{l:1} is satisfied, then the vector field $\xi$ is conformal on the horizontal distribution $\mathcal{H}$. 
\end{theorem}

\begin{proof}
	Since $(M,g,\xi,\lambda)$ is a Ricci soliton and using (\ref{sd:2}) in statement (\ref{p:1}), we get	
	\begin{eqnarray}
	&& \frac{1}{2} (\mathscr{L}_{\xi}g)(X,Y)+Ric^{'}(X^{'},Y^{'})\circ \pi -\frac{1}{2}\{g(\nabla_{X}N,Y)+g(\nabla_{Y}N,X)\} \label{p:2} \\
	&&+2\sum_{i}g(\mathcal{A}_{X}X_{i},\mathcal{A}_{Y}X_{i})+\sum_{j}g(\mathcal{T}_{U_{j}}X,\mathcal{T}_{U_{j}}Y)+\lambda g(X,Y)=0, \notag
	\end{eqnarray}
	where $\{X_i\}$ denotes an orthonormal basis of $\mathcal{H}$, for any $X,Y\in\Gamma(HM)$. Using Lemma \ref{l:1}, the equation (\ref{p:2}) is equivalent to
	\begin{eqnarray}
	\frac{1}{2} (\mathscr{L}_{\xi}g)(X,Y)+Ric^{'}(X^{'},Y^{'})\circ \pi +\lambda g(X,Y)=0. \label{ll:2}
	\end{eqnarray}
	On the other hand, since the Riemannian manifold $(B,g^{'})$ is an Einstein, one can see that $\xi$ is a conformal.	
\end{proof} 

\section{Riemannian Submersions from Ricci Solitons and their Harmonicity}

This section deals with the harmonicity of Riemannian submersion from a Ricci soliton onto a Riemannian manifold. As a tool, we use the tension field and provide a necessary and sufficient condition for which such a submersion $\pi$ is harmonic.

\begin{definition}
	Let $(M,g)$ and $(B,g^{'})$ be $C^{\infty}$-Riemannian manifolds of dimension $m$ and $n$, respectively and $\pi:(M,g)\rightarrow (B,g^{'})$ be a smooth map between Riemannian manifolds. Then $\pi$ is harmonic if and only if the tension field $\tau(\pi)$ of a map $\pi$ vanishes at each point $p\in M$, that is,
	\begin{equation*}
	\tau(\pi)_{p} \ = \ \sum_{i=1}^{m}(\nabla \pi_{*})(e_i,e_i),
	\end{equation*}
	where $\{e_i\}_{1\leq i \leq m}$ is local orthonormal frame around a point $p\in M$ and $\nabla \pi_{*}$ is the second fundamental form of $\pi$, which is defined by
	\begin{eqnarray*}
		\nabla \pi_{*}(E,F) \ \ = \ \ \nabla^{\pi^{-1}TB}_{E}\pi_{*}F-\pi_{*}(\nabla_{E}F)   ,
	\end{eqnarray*}	
	for any vector fields $E,F\in\Gamma(TM)$.	
\end{definition}
For the theory of harmonic maps, we refer to \cite{Eells:1}.  \medskip

Now, we assume that $\pi:(M,g)\rightarrow (B,g^{'})$ is a Riemannian submersion between Riemannian manifolds, such that $(M,g)$ admits a Ricci soliton with the vertical potential field $V$. \medskip

Let $\{e_1,e_2,...,e_m\}$ be an orthonormal basis on $M$, such that $\{e_i\}_{1\leq i \leq r}$ are vertical and $\{e_i\}_{r+1\leq i \leq m}$ are horizontal. Then, it follows that
\begin{eqnarray*}
	\nabla \pi_{*}(e_i,e_i) & = & \sum_{i=1}^m \big( \nabla_{e_i}^{\pi^{-1}TB}\pi_{*}e_i-\pi_{*}(\nabla_{e_i}e_i)  \big) \\
	& = & \sum_{i=1}^r \big( \nabla_{e_i}^{\pi^{-1}TB}\pi_{*}e_i-\pi_{*}(\nabla_{e_i}e_i)\big)+\sum_{i=r+1}^m \big( \nabla_{e_i}^{\pi^{-1}TB}\pi_{*}e_i-\pi_{*}(\nabla_{e_i}e_i)\big) \\
	& = & -\pi_{*}(\sum_{i=1}^r(\nabla_{e_i}e_i)),
\end{eqnarray*}
such that
\begin{equation*}
\nabla_{e_i}^{\pi^{-1}TB}\pi_{*}e_i \ = \ \pi_{*}(\nabla_{e_i}e_i), \ \ \ \ r+1 \leq i \leq m
\end{equation*}
and since $\{e_i\}_{1\leq i \leq r}$ are vertical,
\begin{equation*}
\nabla_{e_i}^{\pi^{-1}TB}\pi_{*}e_i \ = \ 0.
\end{equation*}
Consequently, one has
\begin{eqnarray}
\tau(\pi) & = & \sum_{i=1}^r\nabla \pi_{*}(e_i,e_i) \ = \ -\pi_{*}(\sum_{i=1}^r\mathcal{T}_{e_i}e_i) \label{c:2}\\
& = & -\pi_{*}N. \notag
\end{eqnarray}

Considering Theorem \ref{tt:1} and the statement (\ref{c:2}), we get the theorem as follows:

\begin{theorem}
	Let $(M,g,\xi,\lambda)$ be a Ricci soliton with the horizontal potential field $\xi$ and $\pi:(M,g)\rightarrow (B,g')$ be a Riemannian submersion between Riemannian manifolds. The Riemannian submersion $\pi$ is harmonic, if any two conditions imply the third one::
	
	\noindent (i) The distributions $\mathcal{H}$ or $\mathcal{V}$ are parallel, \\
	\noindent (ii) Any fiber of $\pi$ is minimal,\\
	\noindent (iii) The horizontal distribution $\mathcal{H}$ is integrable and any fiber of $\pi$ has constant scalar curvature $-\lambda r$.	
\end{theorem}

\begin{proof}
	In case the horizontal or vertical distributions are parallel, from Lemma \ref{l:1}, the fundamental tensor field $\mathcal{T}$ vanishes, identically. Similarly, if any fiber of $\pi$ is minimal, from (\ref{bbb:1}), the tensor field $N$ vanishes, identically. Then, considering the equality (\ref{c:2}), the tension field $\tau$ is zero, which gives the Riemannian submersion $\pi$ is harmonic. On the other hand, since $(M,g,\xi,\lambda)$ is a Ricci soliton, using (\ref{p:1}), we get
	\begin{equation}
	\frac{1}{2}(\mathscr{L}_\xi g)(U,W)+Ric(U,W)+\lambda g(U,W)=0, \label{cl}
	\end{equation}
	for any $U,W\in\Gamma(VM)$. By tracing (\ref{cl}), one has
	\begin{eqnarray}
	\sum_{j=1}^{r} \big\{ g(\nabla_{U_{j}}\xi,U_j)+
	Ric(U_j,U_j)+\lambda g(U_j,U_j) \big\} & = & 0,
	\end{eqnarray}
	where $\{U_j\}$ is an orthononormal frame on vertical distribution. By considering the equality (\ref{sd:1}), it follows
	\begin{eqnarray}
	&& \sum_{j=1}^{r} \big\{ g(\nabla_{U_j}\xi,U_j)+\hat{Ric}(U_j,U_j)- \sum_{i=1}^{n} g((\nabla_{X_{i}}\mathcal{T})(U_{j},U_{j}),X_{i})  \label{cv:1} \\
	&& -\sum_{i=1}^{n}g(\mathcal{A}_{X_i}U_j,\mathcal{A}_{X_i}U_j)+\lambda g(U_j,U_j) \big\}+\Vert N \Vert^{2}=0. \notag
	\end{eqnarray}
	Using  the equalities (2), (5) and (11) in (\ref{cv:1}), one has
	\begin{eqnarray}
	&& \sum_{j=1}^{r} \big\{ \hat{Ric}(U_j,U_j)-\check{\delta}(N)-\sum_{i=1}^{n}g(\mathcal{A}_{X_i}U_j,\mathcal{A}_{X_i}U_j)+\lambda g(U_j,U_j)\big\} \label{h} \\  
	&& -g(N,\xi)+\Vert N \Vert^{2}=0. \notag
	\end{eqnarray}
	Since the horizontal distribution $\mathcal{H}$ is integrable and any fiber of $\pi$ has constant scalar curvature $-\lambda r$, from Eq. (\ref{h}) it follows that
	\begin{eqnarray}
	-\check{\delta}(N)+g(N,\xi-N) & = & 0,
	\end{eqnarray}
	Because the conditions (i) or (ii) is satisfied, the fundamental tensor $\mathcal{T}$ vanishes identically. Considering the Eq. (10), the tensor field $N$ vanishes, identically. From (\ref{c:2}), the tension field $\tau$ is zero which is nothing but $\pi$ is harmonic. Therefore, the proof is completed.
\end{proof}\\

\noindent	\textbf{Acknowledgement}s: This work is supported by 1001-Scientific and Technological Research Projects Funding Program of TUBITAK project number 117F434.

\end{document}